\documentclass[reqno,a4paper]{amsart}
\usepackage{amssymb}%
\usepackage[hmarginratio=1:1]{geometry}%
\usepackage{ifpdf}
\usepackage{listings}
\usepackage{color,xcolor}
\usepackage{colortbl}
\usepackage{multirow}
\usepackage{booktabs}
\usepackage{float}
\usepackage{graphicx}
\usepackage{algorithmic}
\usepackage[linesnumbered,ruled,vlined]{algorithm2e}
\usepackage{subfig}
\usepackage{epstopdf}
\usepackage[colorlinks,linkcolor=red,anchorcolor=red,citecolor=red]{hyperref}
\usepackage[numbers,sort&compress]{natbib}
\usepackage{enumitem}

\theoremstyle{plain}
\newtheorem{theorem}{Theorem}[section]

\newtheorem{corollary}{Corollary}[section]

\theoremstyle{remark}
\newtheorem{remark}{Remark}[section]
\theoremstyle{definition}

\numberwithin{equation}{section}
\allowdisplaybreaks[4]

\begin{document}
\title[]{Recurrence relations and applications for the Maclaurin coefficients of squared and cubic hypergeometric functions}
\author{Zhong-Xuan Mao, Jing-Feng Tian*}

\address{Zhong-Xuan Mao \\
Hebei Key Laboratory of Physics and Energy Technology\\
Department of Mathematics and Physics\\
North China Electric Power University \\
Yonghua Street 619, 071003, Baoding, P. R. China}
\email{maozhongxuan000\symbol{64}gmail.com}

\address{Jing-Feng Tian\\
Hebei Key Laboratory of Physics and Energy Technology\\
Department of Mathematics and Physics\\
North China Electric Power University\\
Yonghua Street 619, 071003, Baoding, P. R. China}
\email{tianjf\symbol{64}ncepu.edu.cn}

\begin{abstract}
In this paper, we present and prove that the coefficients $u_n$ and $v_n$ in the series expansions
$F^2(a,b;c;z) = \sum_{n=0}^\infty u_n z^n$ and $F^3(a,b;c;z) = \sum_{n=0}^\infty v_n z^n$
($a,b,c,z \in \mathbb{C}$ and $-c \notin \mathbb{N} \cup \{0\}$) satisfy second- and third-order linear recurrence relations, respectively, where $F(a,b;c;x)$ denotes the Gaussian hypergeometric function and $\mathbb{C}$ is the complex plane.
Our results provide recurrence relations for the Maclaurin coefficients of the squares and cubes of several classical special functions in the complex domain, including zero-balanced Gauss hypergeometric functions, elliptic integrals, as well as classical orthogonal polynomials such as Chebyshev, Legendre, Gegenbauer, and Jacobi polynomials.
As applications, we first establish the monotonicity of a function involving Gauss hypergeometric functions and then present a new proof of the well-known Clausen's formula.
\end{abstract}

\thanks{*Corresponding author: Jing-Feng Tian(tianjf\symbol{64}ncepu.edu.cn)}

\maketitle

\tableofcontents

\section{Introduction}

For any $a,b,c,z\in\mathbb{C}$, the Gaussian hypergeometric function is defined by
\begin{equation} \label{F-ser}
F(a,b;c;z) = \sum_{k=0}^\infty \frac{(a)_k (b)_k}{(c)_k \, k!} z^k, \quad c \notin \mathbb{N} \cup \{0\}, \quad |z| <1.
\end{equation}
As one of the most important special functions, it plays a fundamental role in mathematics and physics, possessing not only theoretical significance, but also widespread applications.

Benefiting from mature tools in series analysis, the series representation \eqref{F-ser} provides a convenient framework for investigating the analytic properties of the Gaussian hypergeometric function, including monotonicity \cite{Wang-IJPAM-2023,Mao-MS-2025}, absolute monotonicity \cite{Wu-BIMS-2024,Zhao-ARXIV-2025}, and inequalities \cite{Tan-MS-2021,Zhao-RACSAM-2022}. Motivated by this idea, in order to treat the square and the cube of the Gaussian hypergeometric function in a similar manner, we also seek their series representations. One possible approach is to employ the Cauchy product formula, by which one can obtain
\begin{equation} \label{CPF-F}
F^2(a,b;c;z) = \sum_{n=0}^\infty \sum_{k=0}^n \frac{(a)_k (a)_{n-k} (b)_k (b)_{n-k}}{(c)_k (c)_{n-k}  k! (n-k)!}  z^n, \quad c \notin \mathbb{N} \cup \{0\},\quad |z| <1.
\end{equation}
and
\begin{equation} \label{CPF-F3}
F^3(a,b;c;z) = \sum_{n=0}^{\infty}
\sum_{i=0}^{n}
\sum_{j=0}^{n-i}
\frac{(a)_i (b)_i}{(c)_i\, i!}\,
\frac{(a)_j (b)_j}{(c)_j\, j!}\,
\frac{(a)_{n-i-j} (b)_{n-i-j}}{(c)_{n-i-j}\, (n-i-j)!}
z^n, \quad c \notin \mathbb{N} \cup \{0\}, \quad |z| <1.
\end{equation}
An alternative approach is to use recurrence relations. Suppose
\begin{equation*}
\begin{aligned}
F^2(a,b;c;z) &= \sum_{n=0}^\infty u_n z^n, \quad |z| <1,\\
F^3(a,b;c;z) &= \sum_{n=0}^\infty v_n z^n, \quad |z| <1.
\end{aligned}
\end{equation*}
By deriving the recurrence relations satisfied by $u_n$ and $v_n$, we can determine the series representations of $F^2$ and $F^3$, thereby providing a solid foundation for analyzing their properties. Some researchers have previously pursued the second approach and obtained some interesting and significant results in the real domain.
In 2018, Yang \cite{Yang-JMAA-2018} established recurrence relations for the Maclaurin series coefficients of $(r^\prime)^p K(r)$ and $(r^\prime)^p E(r)$, where $r \in (0,1)$, $K(r)$ and $E(r)$ denote the complete elliptic integrals of the first and second kinds, respectively.
Later, Chen and Zhao \cite{Chen-RACSAM-2021} studied the recurrence relations for the coefficients of $F^2(a,b;c;x)$, where $a,b,c\in\mathbb{R}$ and $x\in (0,1)$.
More recently, in 2025, Yang \cite{Yang-PAMS-2025} investigated the recurrence relations for the coefficients of $(1-\theta x)^{-q} F(a,b;c;x)$ for $\theta,q,a,b,c\in\mathbb{R}$ and $x\in (0,1)$, provided a new proof of Euler's linear transformation formula, and demonstrated the absolute monotonicity of a specific function.

In view of the importance of this class of methods, we establish and prove the recurrence relations satisfied by $u_n$ and $v_n$ in the complex domain, which are second- and third-order linear recurrence relations, respectively.
Our results provide recurrence relations for the Maclaurin coefficients of the squares and cubes of several classical special functions, including zero-balanced Gauss hypergeometric functions, elliptic integrals, as well as classical orthogonal polynomials such as Chebyshev, Legendre, Gegenbauer, and Jacobi polynomials.
As applications, we first establish the monotonicity of a function involving Gauss hypergeometric functions and then present a new proof of the well-known Clausen's formula.

\section{Recurrence relations for the Maclaurin series coefficients of the square of the Gauss hypergeometric function}

\begin{theorem} \label{thm-1}
Let $a, b, c, z \in\mathbb{C}$ and $-c \notin \mathbb{N} \cup\{0\}$. Then we have
\begin{equation*}
F^2(a,b;c;z) = \sum_{n=0}^\infty u_n z^n, \quad |z| < 1,
\end{equation*}
with $u_0 = 1$, $u_1 = 2 a b / c$ and
\begin{equation*}
u_{n+1} = \alpha_0(n) u_n + \alpha_1(n) u_{n-1}, \quad n\geq1,
\end{equation*}
where
\begin{equation*}
\alpha_0(n) = \frac{2 n^3 + 3(a+b+c-1) n^2 +((a+b)(4c-3)+4ab-c+1)n +2 a b (2 c-1)}{(n+1) (c+n) (2 c+n-1)},
\end{equation*}
and
\begin{equation*}
\alpha_1(n) = -\frac{ (2 a+n-1) (2 b+n-1) (a+b+n-1)}{(n+1) (c+n) (2 c+n-1)}.
\end{equation*}
\end{theorem}

\begin{proof}
Let
\begin{equation*}
\begin{aligned}
R_{-1}(n) &= (n+1) (c+n) (2 c+n-1), \\
R_0(n) &= 2 n^3 + 3(a+b+c-1) n^2 +((a+b)(4c-3)+4ab-c+1)n +2 a b (2 c-1), \\
R_1(n) &= -(2 a+n-1) (2 b+n-1) (a+b+n-1).
\end{aligned}
\end{equation*}
We first show that $R_{-1}(n) u_{n+1} = R_0(n) u_n + R_1(n) u_{n-1}$ holds for all $n \ge 1$, which is equivalent to proving
\begin{equation} \label{sum_r0}
\sum_{n=0}^\infty \Big( R_{-1}(n) u_{n+1} - R_0(n) u_n - R_1(n) u_{n-1} \Big) z^n = 0, \quad |z|<1.
\end{equation}

For convenience, we denote $F = F(a,b;c;z)$. Using the identities
\begin{equation*}
\begin{aligned}
& \sum_{n=0}^\infty u_n z^n = F^2, \\
& \sum_{n=0}^\infty n u_n z^n = z (F^2)^\prime, \\
& \sum_{n=0}^\infty n (n-1) u_n z^n = z^2 (F^2)^{\prime\prime}, \\
& \sum_{n=0}^\infty n (n-1)(n-2) u_n z^n = z^3 (F^2)^{\prime\prime\prime},
\end{aligned}
\end{equation*}
we obtain
\begin{equation*}
\sum_{n=0}^\infty R_{-1}(n) u_{n+1} z^n = z^2 (F^2)^{\prime \prime \prime} + 3c z (F^2)^{\prime \prime} + (2c^2-c) (F^2)^{\prime },
\end{equation*}
\begin{eqnarray*}
\begin{aligned}
\sum_{n=0}^\infty R_{0}(n) u_{n} z^n =& 2 z^3 (F^2)^{\prime \prime \prime} + 3 (a+b+c+1) z^2 (F^2)^{\prime \prime} \\
& \quad +2 (2 a (b+c)+2 b c+c) z (F^2)^{\prime} + 2 a b (2 c-1) F^2,
\end{aligned}
\end{eqnarray*}
and
\begin{eqnarray*}
\begin{aligned}
\sum_{n=0}^\infty R_{1}(n) u_{n-1} z^n = & - z^4 (F^2)^{\prime \prime \prime} -3 (a+b+1) z^3 (F^2)^{\prime \prime} \\
& \quad - \left(2 a^2+a (8 b+3)+2 b^2+3 b+1\right) z^2 (F^2)^{\prime}  -4 a b (a+b) z F^2.
\end{aligned}
\end{eqnarray*}
Furthermore, we have
\begin{eqnarray} \label{sum_r}
\begin{aligned}
& \sum_{n=0}^\infty (R_{-1}(n) u_{n+1} - R_0(n) u_{n} - R_1(n)) z^n \\
= & (z^2 - 2 z^3 + z^4) (F^2)^{\prime \prime \prime} + 3 z (- z (a+b+c+1)+ z^2 (a+b+1)+c) (F^2)^{\prime \prime} \\
& + (z^2 (2 a^2+a (8 b+3)+2 b^2+3 b+1)-2 z (2 a (b+c)+2 b c+c)+2 c^2-c) (F^2)^{\prime} \\
& + 2 a b (2 z (a+b)-2 c+1) F^2.
\end{aligned}
\end{eqnarray}

Since $F$ is a solution of the hypergeometric equation, we have
\begin{equation*}
F^{\prime \prime} = \frac{c-(a+b+1) z}{z(1- z)} F^\prime + \frac{ab}{z(1- z)} F,
\end{equation*}
and differentiating this yields
\begin{eqnarray*}
&& F^{\prime \prime \prime} = \frac{
\left(
\begin{aligned}
& a b \big(z (a+b+3)-c-1\big) F\\
& + \Big( z^2 \big(a^2 + a (b+3) + b^2 + 3b + 2\big) + a b z  - 2 c z (a+b+2) + c^2 + c \Big) F^\prime
\end{aligned}
\right)
}{z^2 (z - 1)^2}.
\end{eqnarray*}
Applying the Leibniz formula, we obtain
\begin{eqnarray} \label{F2''}
\begin{aligned}
(F^2)^{\prime \prime} &= 2 F^\prime F^{\prime} + 2 F F^{\prime\prime} \\
&= \frac{2 a b }{z- z^2} F^2 -\frac{2 (c- z (a+b+1))}{z- z^2} F F^{\prime} + 2 F^{\prime} F^{\prime},
\end{aligned}
\end{eqnarray}
and
\begin{eqnarray} \label{F2'''}
\begin{aligned}
 (F^2)^{\prime \prime \prime}
&= 6 F^\prime F^{\prime \prime} + 2 F F^{\prime\prime\prime} \\
&= \frac{2 a b ( z (a+b+3)-c-1)}{z^2 (z-1)^2} F^2 \\
& \quad + \frac{2 \left( z^2 \left(-2 a b+a (a+3)+b^2+3 b+2\right)+4 a bz-2 c z (a+b+2)+c^2+c\right)}{z^2 ( z-1)^2} F F^{\prime} \\
& \quad -\frac{6 (c- z (a+b+1))}{z- z^2} F^{\prime} F^{\prime}.
\end{aligned}
\end{eqnarray}
Substituting \eqref{F2''} and \eqref{F2'''} into \eqref{sum_r}, we observe that the coefficients of $F^2$, $F F^{\prime}$, and $F^{\prime} F^{\prime}$ all vanish, which proves \eqref{sum_r0}.

Finally, according to the Cauchy product formula \eqref{CPF-F}, we obtain $u_0 = 1$ and $u_1 = 2 a b / c$.
\end{proof}


\begin{remark}
The first few terms of the Maclaurin series for $F^2(a,b;c;z)$ are given by
\begin{equation*}
\begin{aligned}
F^2(a,b;c;z) = &  1 + \frac{2 a b}{c} z + \frac{a b (c (2 a b+a+b+1)+a b)}{c^2 (c+1)} z^2 \\
& \quad + \frac{2 a (a+1) b (b+1) (c (2 a b+a+b+2)+3 a b)}{3 c^2 (c+1) (c+2)} z^3 + \cdots.
\end{aligned}
\end{equation*}
\end{remark}

By setting $c = a + b$ in Theorem~\ref{thm-1}, we obtain the recurrence relation for the Maclaurin series coefficients of the square of the zero-balanced Gaussian hypergeometric function $F(a,b;a+b,z)$.

\begin{corollary} \label{cor1-1}
Let $a,b,z\in\mathbb{C}$ and $-(a+b) \notin \mathbb{N} \cup\{0\}$. Then we have
\begin{equation*}
F^2(a,b;a+b;z) = \sum_{n=0}^\infty u_n z^n, \quad |z|<1,
\end{equation*}
with $u_0 = 1$, $u_1 = 2 a b / (a+b)$ and
\begin{equation*}
u_{n+1} = \alpha_0(n) u_n + \alpha_1(n) u_{n-1}, \quad n\geq1,
\end{equation*}
where
\begin{equation*}
\alpha_0(n) = \frac{(2 a+2 b+2 n-1) (n^2 + (2 a+2 b-1)n +2 a b)}{(n+1) (a+b+n
) (2 (a+ b)+n-1)},
\end{equation*}
and
\begin{equation*}
\alpha_1(n) = -\frac{ (2 a+n-1) (2 b+n-1) (a+b+n-1)}{(n+1) (a+b+n) (2 (a+b)+n-1)}.
\end{equation*}
\end{corollary}

\begin{remark}
The first few terms of the Maclaurin series for $F^2(a,b;a+b;z)$ are given by
\begin{equation*}
\begin{aligned}
F^2(a,b;a+b;z) = &  1 + \frac{2 a b}{a+b} z + \frac{a b \left(a^2 (2 b+1)+a (b+1) (2 b+1)+b^2+b\right)}{(a+b)^2 (a+b+1)} z^2 \\
& \quad + \frac{2 a (a+1) b (b+1) \left(a^2 (2 b+1)+a (b+2) (2 b+1)+b (b+2)\right)}{3 (a+b)^2 (a+b+1) (a+b+2)} z^3 + \cdots.
\end{aligned}
\end{equation*}
\end{remark}

The complete elliptic integral of the first kind, $K(z)$, is defined by
\begin{equation*}
K(z) = \int_0^{\pi/2} \frac{d\theta}{\sqrt{1 - z^2 \sin^2 \theta}}, \quad |z| <1.
\end{equation*}
Its connection with the hypergeometric function $F(a,b;c;z)$ is given by
\begin{equation*}
K(z) = \frac{\pi}{2} F\Bigl(\tfrac{1}{2}, \tfrac{1}{2}; 1; z^2\Bigr).
\end{equation*}
By taking $a = b = \tfrac{1}{2}$ in Corollary~\ref{cor1-1}, we then obtain the recurrence relation for the Maclaurin series coefficients of $K^2(z)$.

\begin{corollary} \label{cor1-1-1}
Let $z\in \mathbb{C}$. Then
\begin{equation*}
K^2(z) = \frac{\pi^2}{4} F^2(\tfrac{1}{2},\tfrac{1}{2};1;z^2) =  \frac{\pi^2}{4} \sum_{n=0}^\infty u_n z^{2n},  \quad |z| < 1,
\end{equation*}
with $u_0 = 1$, $u_1 = 1/2$ and
\begin{equation*}
u_{n+1} = \frac{n (n (2 n+3)+2)+\frac{1}{2}}{(n+1)^3} u_n -\frac{n^3}{(n+1)^3} u_{n-1}, \quad n\geq1.
\end{equation*}
\end{corollary}

\begin{remark}
The first few terms of the Maclaurin series for $K^2(z)$ are given by
\begin{equation*}
K^2(z) = \frac{\pi^2}{4} \Biggl( 1 + \frac{1}{2} z^2 + \frac{11}{32} z^4 + \frac{17}{64} z^6 + \frac{1787}{8192} z^8 + \cdots \Biggr).
\end{equation*}
\end{remark}

The complete elliptic integral of the second kind, $E(z)$, is defined by
\begin{equation*}
E(z) = \int_0^{\pi/2} \sqrt{1 - z^2 \sin^2 \theta} \, d\theta, \quad |z| <1.
\end{equation*}
Its connection with the hypergeometric function $F(a,b;c;z)$ is given by
\begin{equation*}
E(z) = \frac{\pi}{2} F\Bigl(-\tfrac{1}{2}, \tfrac{1}{2}; 1; z^2\Bigr).
\end{equation*}
By taking $a = -\tfrac{1}{2}$, $b = \tfrac{1}{2}$ and $c=1$ in Theorem \ref{thm-1}, we then obtain the recurrence relation for the Maclaurin series coefficients of $E^2(z)$.

\begin{corollary} \label{cor1-1-2}
Let $z\in \mathbb{C}$. Then
\begin{equation*}
E^2(z) = \frac{\pi^2}{4} F^2(-\tfrac{1}{2},\tfrac{1}{2};1;z^2) =  \frac{\pi^2}{4} \sum_{n=0}^\infty u_n z^{2n},  \quad |z| <1,
\end{equation*}
with $u_0 = 1$, $u_1 = -1/2$ and
\begin{equation*}
u_{n+1} = \frac{4 n^3 - 2 n -1}{2(n+1)^3} u_n - \frac{(n-2) (n-1) n}{(n+1)^3} u_{n-1}, \quad n\geq1.
\end{equation*}
\end{corollary}

\begin{remark}
The first few terms of the Maclaurin series for $E^2(z)$ are given by
\begin{equation*}
E^2(z) = \frac{\pi^2}{4} \Biggl( 1 - \frac{1}{2} z^2 - \frac{1}{32} z^4 - \frac{1}{64} z^6 - \frac{77}{8192} z^8 - \cdots \Biggr).
\end{equation*}
\end{remark}

Taking $a=-m$, $b=m$ and $c=\tfrac{1}{2}$ in Theorem \ref{thm-1}, we establish a recurrence relation for the Maclaurin series coefficients of the square of Chebyshev polynomials $T_m(z)$.
\begin{corollary} \label{cor1-2}
Let $m$ be non-negative integer and $z\in \mathbb{C}$. Then
\begin{equation*}
T_m^2(1-2z) = F^2(-m,m;\tfrac1{2};z) = \sum_{n=0}^\infty u_n z^{n},  \quad |z| < 1,
\end{equation*}
with $u_0 = 1$, $u_1 = -4 m^2$ and
\begin{equation*}
u_{n+1} = \frac{-8 m^2+4 n^2-3 n+1}{(2n+1)(n+1)} u_n -\frac{2 (n-1) (-2 m+n-1) (2 m+n-1)}{n (n+1) (2 n+1)} u_{n-1}, \quad n\geq1.
\end{equation*}
\end{corollary}

Taking $a=-m$, $b=m+1$ and $c=1$ in Theorem \ref{thm-1}, we establish a recurrence relation for the Maclaurin series coefficients of the square of Legendre polynomials $P_m(z)$.
\begin{corollary} \label{cor1-3}
Let $m$ be non-negative integer and $z\in \mathbb{C}$. Then
\begin{equation*}
P_m^2(1-2z) = F^2(-m,m+1;1;z) = \sum_{n=0}^\infty u_n z^{n},  \quad |z| < 1,
\end{equation*}
with $u_0 = 1$, $u_1 = -2 m (m+1)$ and
\begin{equation*}
u_{n+1} = \frac{(2 n+1) \left(-2 m (m+1)+n^2+n\right)}{(n+1)^3} u_n - \frac{(2 m+1)^2 n-n^3}{(n+1)^3} u_{n-1}, \quad n\geq1.
\end{equation*}
\end{corollary}

Taking $a=-m$, $b=m+2\alpha$, and $c=\alpha+\tfrac{1}{2}$ in Theorem~\ref{thm-1}, we establish a recurrence relation for the Maclaurin series coefficients of the square of Gegenbauer polynomials $C_m^{(\alpha)}(z)$.

\begin{corollary} \label{cor1-4}
Let $m$ be a non-negative integer, $\alpha>0$ and $z\in \mathbb{C}$. Then
\begin{equation*}
(C_m^{(\alpha)}(1-2z))^2
= \Big(\frac{(2\alpha)_{n}}{n!}\Big)^2 F^2(-m,\, m+2\alpha;\alpha+\tfrac{1}{2};z)
= \Big(\frac{(2\alpha)_{n}}{n!}\Big)^2  \sum_{n=0}^\infty u_n z^{n}, \quad |z| <1,
\end{equation*}
with $u_0 = 1$, $u_1 = -4 m (2 \alpha +m)/(2 \alpha +1)$ and
\begin{equation*}
\begin{aligned}
u_{n+1} =& \frac{n \left(2 \alpha  (8 \alpha -3)-8 m^2-16 \alpha  m+1\right)-8 \alpha  m (2 \alpha +m)+4 n^3+3 (6 \alpha -1) n^2}{(n+1) (2 \alpha +n) (2 \alpha +2 n+1)} u_n \\
& -\frac{2 (-2 m+n-1) (2 \alpha +n-1) (4 \alpha +2 m+n-1)}{(n+1) (2 \alpha +n) (2 \alpha +2 n+1)} u_{n-1},
\quad n\ge 1.
\end{aligned}
\end{equation*}
\end{corollary}

Taking $a=-m$, $b=m+\alpha+\beta+1$, and $c=\alpha+1$ in Theorem~\ref{thm-1}, we establish a recurrence relation for the Maclaurin series coefficients of the square of Jacobi polynomials $P_m^{(\alpha,\beta)}(z)$.

\begin{corollary} \label{cor1-5}
Let $m$ be a non-negative integer, $\alpha,\beta>0$ and $z\in \mathbb{C}$. Then
\begin{equation*}
\begin{aligned}
(P_m^{(\alpha,\beta)}(1-2z))^2
& = \Big(\frac{(\alpha+1)_{n}}{n!}\Big)^2 F^2(-m,m+\alpha+\beta+1;\alpha+1;z) \\
& = \Big(\frac{(\alpha+1)_{n}}{n!}\Big)^2  \sum_{n=0}^\infty u_n z^{n}, \quad |z| < 1,
\end{aligned}
\end{equation*}
with $u_0 = 1$, $u_1 = -2 m ( m + \alpha +\beta+1)/(\alpha +1)$ and
\begin{equation*}
\begin{aligned}
u_{n+1} =& \frac{\left(
\begin{aligned}
& -2 m^2 (2 \alpha +2 n+1)-2 m (\alpha +\beta +1) (2 \alpha +2 n+1) \\
& \quad +n \left(4 \alpha  (\alpha +\beta +1)+\beta +2 n^2+3 n (2 \alpha +\beta +1)+1\right)
\end{aligned}
\right)}{(n+1) (\alpha +n+1) (2 \alpha +n+1)} u_n \\
& + \frac{(2 m-n+1) (\alpha +\beta +n) (2 \alpha +2 \beta +2 m+n+1)}{(n+1) (\alpha +n+1) (2 \alpha +n+1)} u_{n-1},
\quad n\ge 1.
\end{aligned}
\end{equation*}
\end{corollary}

Taking $a=\tfrac{1+\alpha}{2}$, $b=\tfrac{1-\alpha}{2}$ and $c=\tfrac{3}{2}$ in Theorem \ref{thm-1}, we have the following corollary.

\begin{corollary} \label{cor1-6}
Let $\alpha, z \in \mathbb{C}$. Then
\begin{equation*}
\sin^2(\alpha \arcsin z) = \alpha^2 z^2 F^2(\tfrac{1+\alpha}{2},\tfrac{1-\alpha}{2};\tfrac{3}{2};z^2) = \alpha^2 \sum_{n=0}^\infty u_n z^{2n+2}, \quad |z| < 1,
\end{equation*}
with $u_0 = 1$, $u_1 = (1-\alpha^2)/3$ and
\begin{equation*}
u_{n+1} = \frac{-2 \alpha ^2+4 n^2+5 n+2}{2 n^2+7 n+6} u_n -\frac{2 n (n-\alpha ) (\alpha +n)}{(n+1) (n+2) (2 n+3)} u_{n-1}, \quad n\ge 1.
\end{equation*}
\end{corollary}

Taking $a=\tfrac{\alpha}{2}$, $b=-\tfrac{\alpha}{2}$ and $c=\tfrac{1}{2}$ in Theorem \ref{thm-1}, we have the following corollary.

\begin{corollary} \label{cor1-7}
Let $\alpha,z \in \mathbb{C}$. Then
\begin{equation*}
\cos^2(\alpha \arcsin z) = F^2(\tfrac{\alpha}{2},-\tfrac{\alpha}{2};\tfrac{1}{2};z^2) = \sum_{n=0}^\infty u_n z^{2n}, \quad |z| < 1,
\end{equation*}
with $u_0 = 1$, $u_1 = -\alpha^2$ and
\begin{equation*}
u_{n+1} = \frac{-2 \alpha ^2+4 n^2-3 n+1}{2 n^2+3 n+1} u_n -\frac{2 (n-1) (-\alpha +n-1) (\alpha +n-1)}{n (n+1) (2 n+1)} u_{n-1}, \quad n\ge 1.
\end{equation*}
\end{corollary}

For any integers $n_1,n_2,n_3$, the function
$F(a+n_1,b+n_2;c+n_3;z)$ is called a contiguous function of $F(a,b;c;z)$.
Gauss proved that any three contiguous hypergeometric functions are linearly related,
giving rise to a total of 15 independent linear relations, known as Gauss's contiguous relations.
From Theorem \ref{thm-1}, we provide the power series representation of $F(a,b;c;z)$ along with that of one of its contiguous functions.

\begin{corollary}
Let $a,b,c,z\in\mathbb{C}$ and $-c \notin \mathbb{N} \cup\{0\}$. Then we have
\begin{equation*}
F(a,b;c;z) F(a+1,b+1;c+1;z)  = \sum_{n=0}^\infty u_n z^n, \quad |z| <1,
\end{equation*}
with $u_0 = 1$, $u_1 = \tfrac{a b}{c}+\tfrac{(a+1) (b+1)}{c+1}$, and
\begin{equation*}
\frac{u_{n}}{n+1} = \alpha_0(n) \frac{u_{n-1}}{n} + \alpha_1(n) \frac{u_{n-2}}{n-1}, \quad n\geq2,
\end{equation*}
where $\alpha_0(n)$ and $\alpha_1(n)$ are defined in Theorem \ref{thm-1}.
\end{corollary}

\section{Recurrence relations for the Maclaurin series coefficients of the cube of the Gauss hypergeometric function}

\begin{theorem} \label{thm-2}
Let $a, b, c, z \in\mathbb{C}$ and $-c \notin \mathbb{N} \cup\{0\}$. Then we have
\begin{equation*}
F^3(a,b;c;z) = \sum_{n=0}^\infty v_n z^n, \quad |z| < 1,
\end{equation*}
with $v_0 = 1$, $v_1 = 3 a b / c$, $v_2 = \frac{3 a^2 b^2}{c^2}+\frac{3 a (a+1) (b+1) b}{2 c (c+1)}$ and
\begin{equation*}
v_{n+1} = \beta_0(n) v_n + \beta_1(n) v_{n-1} + \beta_2(n) v_{n-2}, \quad n\geq1,
\end{equation*}
where
\begin{equation} \label{beta}
\begin{aligned}
\beta_0(n) &= \frac{\left(
\begin{aligned}
& 3 n^4 + 6 (a+b+2 c-2) n^3 + (2 a (5 b+11 c-9)+11 c (2 b+c)-18 b-26 c+15) n^2  \\
& \quad +\left(c^2 (18 a+18 b-7)+c (a (30 b-29)-29 b+13)+4 a (3-5 b)+12 b-6\right) n \\
& \quad  + 3 a b (c (6 c-7)+2)
\end{aligned} \right)}{(n+1) (c+n) (2 c+n-1) (3 c+n-2)}, \\
\beta_1(n) &= -\frac{\left(
\begin{aligned}
& 3 n^4 + 6 (2 a+2 b+c-3) n^3 \\
& \quad + \left(11 a^2+42 a b+22 a c-54 a+11 b^2+22 b c-54 b-22 c+40\right)n^2 \\
& \quad +3 \left(10 a^2 b+6 a^2 c-11 a^2+10 a b^2+22 a b c-42 a b-17 a c+26 a \right)n\\
& \quad +3 \left( 6 b^2 c-11 b^2-17 b c+26 b+9 c-13\right)n \\
& \quad + 9 a^2 b^2+36 a^2 b c-45 a^2 b-18 a^2 c+22 a^2+36 a b^2 c-45 a b^2-72 a b c \\
& \quad +87 a b+29 a c-36 a-18 b^2 c+22 b^2+29 b c-36 b-11 c+14
\end{aligned}
\right)}{(n+1) (c+n) (2 c+n-1) (3 c+n-2)}, \\
\beta_2(n) &= \frac{(3 a+n-2) (3 b+n-2) (2 a+b+n-2) (a+2 b+n-2)}{(n+1) (c+n) (2 c+n-1) (3 c+n-2)}.
\end{aligned}
\end{equation}
\end{theorem}

\begin{proof}
The proof is similar to that of Theorem~\ref{thm-1}, and is therefore omitted.
\end{proof}

By taking $a = b = \tfrac{1}{2}$ and $c=1$ in Theorem \ref{thm-1}, we then obtain the recurrence relation for the Maclaurin series coefficients of the cubic power of the complete elliptic integral of the first kind $K(z)$, $z\in \mathbb{C}$.

\begin{corollary} \label{cor1-1-1}
Let $z\in \mathbb{C}$. Then
\begin{equation*}
K^3(z) = \frac{\pi^3}{8} F^3(\tfrac{1}{2},\tfrac{1}{2};1;z^2) =  \frac{\pi^3}{8} \sum_{n=0}^\infty v_n z^{3n},  \quad |z| < 1,
\end{equation*}
with $v_0 = 1$, $v_1 = 3/4$, $v_2 = 39/64$ and
\begin{equation*}
v_{n+1} = \frac{2 n (n+1) (6 n (n+1)+7)+3}{4 (n+1)^4} v_n - \frac{48 n^4+32 n^2+1}{16 (n+1)^4} v_{n-1} + \frac{(1-2 n)^4}{16 (n+1)^4} v_{n-2}, \quad n \geq 2.
\end{equation*}
\end{corollary}

By taking $a = -\tfrac{1}{2}$, $b = \tfrac{1}{2}$ and $c=1$ in Theorem \ref{thm-1}, we then obtain the recurrence relation for the Maclaurin series coefficients of the cubic power of the complete elliptic integral of the second kind $E(z)$, $z\in \mathbb{C}$.

\begin{corollary} \label{cor1-2-1}
Let $z\in \mathbb{C}$. Then
\begin{equation*}
E^3(z) = \frac{\pi^3}{8} F^3(-\tfrac{1}{2},\tfrac{1}{2};1;z^2) =  \frac{\pi^3}{8} \sum_{n=0}^\infty v_n z^{3n},  \quad |z| < 1,
\end{equation*}
with $v_0 = 1$, $v_1 = -3/4$, $v_2 = 3/64$ and
\begin{equation*}
\begin{aligned}
v_{n+1} = & \frac{2 n \left(6 n^3-5 n-5\right)-3}{4 (n+1)^4} v_n -\frac{8 n (2 n (3 (n-4) n+13)-9)+29}{16 (n+1)^4} v_{n-1} \\
& \quad + \frac{(2 n-7) (2 n-5) (2 n-3) (2 n-1)}{16 (n+1)^4} v_{n-2}, \quad n \geq 2.
\end{aligned}
\end{equation*}
\end{corollary}

Taking $a=-m$, $b=m$ and $c=\tfrac{1}{2}$ in Theorem \ref{thm-2}, we establish a recurrence relation for the Maclaurin series coefficients of the cubic of Chebyshev polynomials $T_m(z)$.
\begin{corollary} \label{cor2-2}
Let $m$ be non-negative integer and $z\in \mathbb{C}$. Then
\begin{equation*}
T_m^3(1-2z) = F^3(-m,m;\tfrac1{2};z) = \sum_{n=0}^\infty v_n z^{n},  \quad |z| < 1,
\end{equation*}
with $v_0 = 1$, $v_1 = - 6 m^2$, $v_2 = 2 (7 m^4-m^2)$ and
\begin{equation*}
\begin{aligned}
v_{n+1} =& \frac{-20 m^2+6 n^2-9 n+5}{2 n^2+3 n+1} v_n \\
&  + \frac{-36 m^4+20 m^2 (n-1) (4 n-5)-2 (n-1) (2 n (3 (n-4) n+17)-17)}{n (n+1) (2 n-1) (2 n+1)} v_{n-1} \\
& + \frac{4 \left(9 m^4-10 m^2 (n-2)^2+(n-2)^4\right)}{n (n+1) (2 n-1) (2 n+1)} v_{n-2}, \quad n\geq1.
\end{aligned}
\end{equation*}
\end{corollary}

Taking $a=-m$, $b=m+1$ and $c=1$ in Theorem \ref{thm-2}, we establish a recurrence relation for the Maclaurin series coefficients of the cubic of Legendre polynomials $P_m(z)$, $z\in \mathbb{C}$.

\begin{corollary} \label{cor2-3}
Let $m$ be non-negative integer and $z\in \mathbb{C}$. Then
\begin{equation*}
P_m^3(1-2z) = F^3(-m,m+1;1;z) = \sum_{n=0}^\infty v_n z^{n},  \quad |z| < 1,
\end{equation*}
with $v_0 = 1$, $v_1 = - 3 m(m+1)$, $v_2 = \frac{3}{4} (5 m^4+10 m^3+3 m^2-2 m)$ and
\begin{equation*}
\begin{aligned}
& v_{n+1} \\
=& \frac{3 (n+1) (2 n-1)-20 m (m+1)}{(n+1) (2 n+1)} v_n \\
&  + \frac{4 \left((20 m (m+1)+3) n^2-3 (5 m (m+1)+1) n-m (m+1) (3 m+1) (3 m+2)-3 n^4+3 n^3\right)}{n (n+1) (2 n-1) (2 n+1)} v_{n-1} \\
& + \frac{4 (m-n+1) (3 m-n+2) (m+n) (3 m+n+1)}{n (n+1) (2 n-1) (2 n+1)} v_{n-2}, \quad n\geq1.
\end{aligned}
\end{equation*}
\end{corollary}

Taking $a=-m$, $b=m+2\alpha$ and $c=\alpha+\tfrac{1}{2}$ in Theorem~\ref{thm-1}, we establish a recurrence relation for the Maclaurin series coefficients of the cubic of Gegenbauer polynomials $C_m^{(\alpha)}(z)$, $z\in \mathbb{C}$.

\begin{corollary} \label{cor1-4}
Let $m$ be a non-negative integer, $\alpha>0$ and $z\in \mathbb{C}$. Then
\begin{equation*}
(C_m^{(\alpha)}(1-2z))^3
= \Big(\frac{(2\alpha)_{n}}{n!}\Big)^3 F^3(-m,\, m+2\alpha;\alpha+\tfrac{1}{2};z)
= \Big(\frac{(2\alpha)_{n}}{n!}\Big)^3  \sum_{n=0}^\infty v_n z^{n}, \quad |z| <1,
\end{equation*}
with $v_0 = 1$, $v_1 = -5 m (2 \alpha +m)/(2 \alpha +1)$, $v_2 = \frac{3 m^2 (2 \alpha +m)^2}{\left(\alpha +\frac{1}{2}\right)^2}+\frac{6 (m-1) m (2 \alpha +m+1) (2 \alpha +m)}{(2 \alpha +1) (2 \alpha +3)}$ and
\begin{equation*}
\begin{aligned}
v_{n+1} =& \frac{
\left(
\begin{aligned}
& 12 n^4 + 24 (4 \alpha -1) n^3 + (4 \alpha  (55 \alpha -29)-40 m^2-80 \alpha  m+19 )n^2 \\
& \quad + \left((4 \alpha -1) (4 \alpha  (9 \alpha -5)+5)+20 (1-6 \alpha ) m^2+40 \alpha  (1-6 \alpha ) m\right) n \\
& \quad - 12 \alpha  (6 \alpha -1) m (2 \alpha +m)
\end{aligned}
\right)
}{(n+1) (2 \alpha +n) (2 \alpha +2 n+1) (6 \alpha +2 n-1)} v_n \\
& + \frac{\left(
\begin{aligned}
& -12 n^4 + (60-120 \alpha) n^3 + 4(-88 \alpha ^2+108 \alpha +20 m^2+40 \alpha  m-29)n^2  \\
& \quad 6 (2 \alpha -1)(-24 \alpha ^2+54 \alpha +30 m^2+60 \alpha  m-17) n \\
& \quad + 288 \alpha ^3-440 \alpha ^2+216 \alpha -36 m^4-144 \alpha  m^3+144 \alpha ^2 m^2-360 \alpha  m^2 \\
& \quad +100 m^2+576 \alpha ^3 m-720 \alpha ^2 m+200 \alpha  m-34
\end{aligned}
\right)
}{(n+1) (2 \alpha +n) (2 \alpha +2 n+1) (6 \alpha +2 n-1)} v_{n-1}\\
& + \frac{4 (3 m-n+2) (-2 \alpha +m-n+2) (4 \alpha +m+n-2) (6 \alpha +3 m+n-2)}{(n+1) (2 \alpha +n) (2 \alpha +2 n+1) (6 \alpha +2 n-1)} v_{n-2},
\quad n\ge 2.
\end{aligned}
\end{equation*}
\end{corollary}

Taking $a=-m$, $b=m+\alpha+\beta+1$ and $c=\alpha+1$ in Theorem~\ref{thm-1}, we establish a recurrence relation for the Maclaurin series coefficients of the cubic of Jacobi polynomials $P_m^{(\alpha,\beta)}(z)$, $z\in \mathbb{C}$.

\begin{corollary} \label{cor1-5}
Let $m$ be a non-negative integer, $\alpha,\beta>0$ and $z\in \mathbb{C}$. Then
\begin{equation*}
\begin{aligned}
(P_m^{(\alpha,\beta)}(1-2z))^3
& = \Big(\frac{(\alpha+1)_{n}}{n!}\Big)^3 F^3(-m,m+\alpha+\beta+1;\alpha+1;z) \\
& = \Big(\frac{(\alpha+1)_{n}}{n!}\Big)^3  \sum_{n=0}^\infty v_n z^{n}, \quad |z| < 1,
\end{aligned}
\end{equation*}
with $v_0 = 1$, $u_1 = - 3 m ( m + \alpha +\beta+1)/(\alpha +1)$, $u_2 = \frac{3 m^2 (\alpha +\beta +m+1)^2}{(\alpha +1)^2}+\frac{3 (m-1) m (\alpha +\beta +m+2) (\alpha +\beta +m+1)}{2 (\alpha +1) (\alpha +2)}$ and
\begin{equation*}
\begin{aligned}
& v_{n+1} \\
=& \frac{
\left(
\begin{aligned}
& 3 n^4 + 6 (3 \alpha +\beta +1)n^3 + \\
& \quad +\left(11 \alpha  (3 \alpha +2 \beta +2)+4 (\beta +1)-10 m^2-10 m (\alpha +\beta +1)\right) n^2 \\
& \quad + \left(18 \alpha ^3+18 \alpha ^2 (\beta +1)+7 \alpha  (\beta +1)+\beta\right) n \\
& \quad - \left( 10 (3 \alpha +1) m^2 + 10 (3 \alpha +1) m (\alpha +\beta +1)-1\right) n \\
&  \quad -3 \left(6 \alpha ^2+5 \alpha +1\right) m (\alpha +\beta +m+1)
\end{aligned}
\right)
}{(n+1) (\alpha +n+1) (2 \alpha +n+1) (3 \alpha +n+1)} v_n \\
& + \frac{
\left(
\begin{aligned}
& -3 n^4 - 6 (3 \alpha +2 \beta )n^3 + \\
& \quad + (-33 \alpha ^2+\alpha  (10-44 \beta )-11 \beta ^2+10 \beta +20 m^2+20 m (\alpha +\beta +1)+3) n^2\\
& \quad 3 (-6 \alpha ^3+\alpha (-6 \beta ^2+15 \beta +20 m^2+10 (3 \beta +2) m+3)) n \\
& \quad + 3( \beta  \left(5 \beta +10 m^2+10 (\beta +1) m+1\right)+2 \alpha ^2 (-6 \beta +10 m+5)) n \\
& \quad + (\alpha +\beta ) (3 \alpha  (6 \alpha +6 \beta +1)-4 \beta -1)-9 m^4-18 m^3 (\alpha +\beta +1) \\
& \quad +m^2 (27 \alpha ^2+9 \alpha  (2 \beta -3)-9 \beta  (\beta +3)-11 ) \\
& \quad +m (\alpha +\beta +1) (9 \alpha  (4 \alpha +4 \beta -1)-9 \beta -2)
\end{aligned}
\right)
}{(n+1) (\alpha +n+1) (2 \alpha +n+1) (3 \alpha +n+1)} v_{n-1}\\
& + \frac{(-3 m+n-2) (\alpha +\beta -m+n-2) (3 \alpha +3 \beta +3 m+n-2) (2 (\alpha +\beta -1)+m+n)}{(n+1) (\alpha +n+1) (2 \alpha +n+1) (3 \alpha +n+1)} v_{n-2},
\quad n\ge 2.
\end{aligned}
\end{equation*}
\end{corollary}

Taking $a=\tfrac{1+\alpha}{2}$, $b=\tfrac{1-\alpha}{2}$ and $c=\tfrac{3}{2}$ in Theorem \ref{thm-1}, we have the following corollary.

\begin{corollary} \label{cor2-6}
Let $\alpha,z \in \mathbb{C}$. Then
\begin{equation*}
\sin^3(\alpha \arcsin z) = \alpha^3 z^3 F^3(\tfrac{1+\alpha}{2},\tfrac{1-\alpha}{2};\tfrac{3}{2};z^2) = \alpha^3 \sum_{n=0}^\infty v_n z^{2n+3}, \quad |z| < 1,
\end{equation*}
with $v_0 = 0$, $v_1 = (1-\alpha^2)/2$, $v_2 = (13 \alpha ^4-50 \alpha ^2+37)/120$ and
\begin{equation*}
\begin{aligned}
v_{n+1} = & \frac{-5 \alpha ^2+6 n^2+9 n+5}{2 n^2+9 n+10} v_n \\
& + \frac{-9 \alpha ^4+10 \alpha ^2+20 \alpha ^2 n (4 n+3)-4 n (4 n (3 n (n+1)+2)+3)-1}{4 (n+1) (n+2) (2 n+3) (2 n+5)} v_{n-1} \\
& + \frac{9 \alpha ^4-10 \alpha ^2 (1-2 n)^2+(1-2 n)^4}{4 (n+1) (n+2) (2 n+3) (2 n+5)} v_{n-2}, \quad n\geq 2.
\end{aligned}
\end{equation*}
\end{corollary}

Taking $a=\tfrac{\alpha}{2}$, $b=-\tfrac{\alpha}{2}$ and $c=\tfrac{1}{2}$ in Theorem \ref{thm-1}, we have the following corollary.

\begin{corollary} \label{cor2-7}
Let $\alpha,z \in \mathbb{C}$. Then
\begin{equation*}
\cos^3(\alpha \arcsin z) = F^3(\tfrac{\alpha}{2},-\tfrac{\alpha}{2};\tfrac{1}{2};z^2) = \sum_{n=0}^\infty v_n z^{2n}, \quad |z| < 1,
\end{equation*}
with $v_0 = 0$, $v_1 = -3\alpha^2/2$, $v_2 = (7 \alpha ^4-4 \alpha ^2)/8$ and
\begin{equation*}
\begin{aligned}
v_{n+1} = & \frac{-5 \alpha ^2+6 n^2-9 n+5}{2 n^2+3 n+1} v_n \\
& + \frac{-9 \alpha ^4+20 \alpha ^2 (n-1) (4 n-5)-8 (n-1) (2 n (3 (n-4) n+17)-17)}{4 n (n+1) (2 n-1) (2 n+1)} v_{n-1} \\
& + \frac{9 \alpha ^4-40 \alpha ^2 (n-2)^2+16 (n-2)^4}{4 n (n+1) (2 n-1) (2 n+1)} v_{n-2}, \quad n\geq 2.
\end{aligned}
\end{equation*}
\end{corollary}

\begin{corollary}
Let $a,b,c,z\in\mathbb{C}$ and $-c \notin \mathbb{N} \cup\{0\}$. Then we have
\begin{equation*}
F^2(a,b;c;z) F(a+1,b+1;c+1;z)  = \sum_{n=0}^\infty v_n z^n, \quad |z| <1,
\end{equation*}
with $v_0 = 1$, $v_1 = \frac{2 a b}{c}+\frac{(a+1) (b+1)}{c+1}$, $v_2 = \frac{a^2 b^2}{c^2}+\frac{3 a (a+1) (b+1) b}{c (c+1)}+\frac{(a+1) (a+2) (b+1) (b+2)}{2 (c+1) (c+2)}$ and
\begin{equation*}
\frac{v_{n}}{n+1} = \beta_0(n) \frac{v_{n-1}}{n} + \beta_1(n) \frac{v_{n-2}}{n-1} + \beta_2(n) \frac{v_{n-3}}{n-2}, \quad n\geq3,
\end{equation*}
where $\beta_0(n)$, $\beta_1(n)$ and $\beta_2(n)$ are defined in Theorem \ref{thm-2}.
\end{corollary}

%
%
%

\section{Applications}

\subsection{Monotonicity of a function involving the Gauss hypergeometric function}
\begin{theorem} \label{thm-3}
Let $a,b,c>0$ and $(c-a)(c-b) \geq 0$. Then the function
\begin{equation} \label{F2F}
x \mapsto \frac{F^2(a,b;c;x)}{F(2 a, 2 b; 2c;x)}
\end{equation}
is increasing on $(0,1)$. Moreover, we have $F^2(a,b;c;x) \geq F(2 a, 2 b; 2c;x)$ for all $x\in(0,1)$.
\end{theorem}

\begin{proof}
Let the sequence $\{u_n\}_{n\geq0}$ be determined by the following initial values and recurrence relation: $u_0=1$, $u_1=2ab/c$, and
\begin{equation*}
u_{n+1} = \alpha_0(n) u_n + \alpha_1(n) u_{n-1},
\end{equation*}
where $\alpha_0(n)$ and $\alpha_1(n)$ are given in Theorem \ref{cor1-1}, and let the sequence
$\{v_n\}_{n\geq0}$ be defined by
\begin{equation*}
v_n = \frac{(2 a)_n (2 b)_n}{(2c)_n \, n!}.
\end{equation*}
Then we have
\begin{equation*}
\frac{F^2(a,b;c;x)}{F(2 a, 2 b; 2c;x)} = \frac{\sum_{n=0}^\infty u_n x^n}{\sum_{n=0}^\infty v_n x^n}.
\end{equation*}

We will prove that the sequence $\{u_n/v_n\}_{n\geq0}$ is increasing. Since
\begin{equation*}
\frac{u_1}{v_1} - \frac{u_0}{v_0} = 0 \geq 0,
\end{equation*}
we assume that $\{u_k/v_k\}$ is increasing on $k=n-1$, that is, $\tfrac{u_n}{v_n} \geq \tfrac{u_{n-1}}{v_{n-1}}$. By mathematical induction, it suffices to show that
\begin{equation*}
\frac{u_{n+1}}{v_{n+1}} \geq \frac{u_n}{v_n}.
\end{equation*}
Since $u_n$ and $v_n$ are nonnegative if $a,b,c>0$, this is equivalent to proving
\begin{equation*}
\frac{\alpha_0(n) u_n + \alpha_1(n) u_{n-1}}{u_n} \geq \frac{v_{n+1}}{v_n}.
\end{equation*}
Indeed,
\begin{equation*}
\begin{aligned}
 \frac{\alpha_0(n) u_n + \alpha_1(n) u_{n-1}}{u_n} - \frac{v_{n+1}}{v_n}
& \geq \alpha_0(n) + \alpha_1(n) \frac{v_{n-1}}{v_n} - \frac{v_{n+1}}{v_n} \\
& = \frac{2 n (a-c) (b-c)}{(n+1) (c+n) (2 c+n-1) (2 c+n)} \geq 0.
\end{aligned}
\end{equation*}

According to \cite[Lemma 1]{Biernacki-AUMCS-1955}, we obtain that the function \eqref{F2F} is increasing.
\end{proof}

From Theorem \ref{thm-3}, we have the following corollary.
\begin{corollary}
Let $a,b>0$ and $c > \max\{a, b, 0\}$. Then we have
\begin{equation*}
1< \frac{F^2(a,b;c;x)}{F(2 a, 2 b; 2c;x)} < \frac{\Gamma^2(c)\, \Gamma^2(c-a-b)\, \Gamma(2c-2a)\, \Gamma(2c-2b)}
     {\Gamma^2(c-a)\, \Gamma^2(c-b)\, \Gamma(2c)\, \Gamma(2c-2a-2b)}, \quad x\in(0,1).
\end{equation*}
\end{corollary}

\begin{proof}
From Theorem \ref{thm-3}, we know that the function \eqref{F2F} is increasing on $(0,1)$ with
\begin{equation*}
\lim_{x\to 0^+} \frac{F^2(a,b;c;x)}{F(2 a, 2 b; 2c;x)} = 1,
\end{equation*}
and
\begin{equation*}
\lim_{x\to 1^-} \frac{F^2(a,b;c;x)}{F(2 a, 2 b; 2c;x)} = \frac{\Gamma^2(c)\, \Gamma^2(c-a-b)\, \Gamma(2c-2a)\, \Gamma(2c-2b)}
     {\Gamma^2(c-a)\, \Gamma^2(c-b)\, \Gamma(2c)\, \Gamma(2c-2a-2b)}, \quad x\in(0,1).
\end{equation*}
\end{proof}

Taking $c=a+b$ in Theorem \ref{thm-3}, we have the following corollary.
\begin{corollary} \label{cor3-1}
Let $a,b>0$ and $a+b>\tfrac{1}{2}$. Then the function
\begin{equation} \label{F2F}
x \mapsto \frac{F^2(a,b;a+b;x)}{F(2 a, 2 b; 2(a+b);x)}
\end{equation}
is increasing from $(0,1)$ onto $(1,\infty)$.
\end{corollary}

Taking $a=b=\frac{1}{2}$ in Corollary \ref{cor3-1}, we have the following corollary.

\begin{corollary}
It follows that
\begin{equation*}
x\mapsto \frac{x K^2(x)}{\log (1-x)}
\end{equation*}
is decreasing from $(0,1)$ onto $(-\infty,-\frac{\pi ^2}{4})$.
\end{corollary}

\subsection{A new proof for the Clausen formula}
Clausen's formula \cite{Clausen-JM-1828} is given by
\begin{equation*}
F^2(a,b;a+b+\tfrac{1}{2}) = {_3F_2}(2a,2b,a+b;a+b+\tfrac{1}{2},2a+2b;x),
\end{equation*}
for which we provide a new proof below.

\begin{proof}
Let
\begin{equation*}
w_n = \frac{(2 a)_n (2 b)_n (a+b)_n}{n! \left(a+b+\frac{1}{2}\right)_n (2 a+2 b)_n}.
\end{equation*}
Then
\begin{equation*}
{_3F_2}(2a,2b,a+b;a+b+\tfrac{1}{2},2a+2b;x) = \sum_{n=0}^\infty w_n x^n.
\end{equation*}
Observe that $w_0=1=u_0$ and $w_1=2 a b / (a+b+\tfrac{1}{2})=u_1$.
To show that $w_n = u_n$, it suffices to prove that $w_n$ also satisfies the recurrence relation
\begin{equation}\label{w-r}
w_{n+1} = \alpha_0(n) w_n + \alpha_1(n) w_{n-1}.
\end{equation}

Indeed, noting that
\begin{equation*}
\frac{w_{n+1}}{w_n} = \frac{2 (2 a+n) (2 b+n) (a+b+n)}{(n+1) (2 (a+b)+n) (2 a+2 b+2 n+1)},
\end{equation*}
we have
\begin{equation*}
\begin{aligned}
& w_{n+1} - \alpha_0(n) w_n - \alpha_1(n) w_{n-1} \\
=& \frac{4 (2 a+n-1) (2 a+n) (2 b+n-1) (2 b+n) (a+b+n-1) (a+b+n)}{n (n+1) (2 a+2 b+n-1) (2 (a+b)+n) (2 a+2 b+2 n-1) (2 a+2 b+2 n+1)} w_{n-1} \\
& - \alpha_0(n) \frac{2 (2 a+n-1) (2 b+n-1) (a+b+n-1)}{n (2 a+2 b+n-1) (2 a+2 b+2 n-1)} w_{n-1} - \alpha_1(n) w_{n-1} \\
=& 0.
\end{aligned}
\end{equation*}
\end{proof}

The following identity is known as the Ramanujan--Preece formula.
It was originally stated by Ramanujan without proof and later rigorously proved by Preece
\begin{equation*}
M(a,c;x)\, M(a,c;-x) = {_2F_3}\!\left(a,\, c-a;\, c,\, \tfrac{c}{2},\, \tfrac{c+1}{2};\, \tfrac{x^2}{4}\right).
\end{equation*}

\begin{proof}
Following an argument similar to that of Theorem~\ref{thm-2}, we obtain
\begin{equation*}
M(a,c;x)\, M(a,c;-x) = \sum_{n=0}^\infty u_n x^n,
\end{equation*}
where the coefficients $u_n$ satisfy $u_0=1$, $u_1=0$ and
\begin{equation*}
u_{n+1} = -\frac{(2 a+n-1) (2 a-2 c-n+1)}{(n+1) (c+n-1) (c+n) (2 c+n-1)} \, u_{n-1},
\end{equation*}
which implies
\begin{equation*}
u_{2n} = -\frac{(a+n-1) (a-c-n+1)}{n (c+n-1) (c+2 n-2) (c+2 n-1)} \, u_{2n-2}.
\end{equation*}

On the other hand, noting that
\begin{equation*}
{_2F_3}\!\left(a,\, c-a;\, c,\, \tfrac{c}{2},\, \tfrac{c+1}{2};\, \tfrac{x^2}{4}\right) = \sum_{n=0}^\infty w_{2n} x^{2n} = \sum_{n=0}^\infty \frac{(a)_n (c-a)_n}{(c)_n \left(\frac{c}{2}\right)_n \left(\frac{c+1}{2}\right)_n} \frac{1}{4^n n!} x^{2n},
\end{equation*}
where $w_n$ is defined by
\begin{equation*}
w_{2n} = \frac{(a)_n (c-a)_n}{(c)_n \left(\frac{c}{2}\right)_n \left(\frac{c+1}{2}\right)_n} \frac{1}{4^n n!}
\end{equation*}
satisfying $w_0=1$ and
\begin{equation*}
w_{2n} = \frac{(a+n-1) (-a+c+n-1)}{n (c+n-1) (c+2 n-2) (c+2 n-1)} \, w_{2n-2}.
\end{equation*}

Thus, $u_n = w_n$, which is the desired conclusion.
\end{proof}

\end{document}